\documentclass[12pt]{article}
\usepackage{amsmath, amsfonts, amsthm, amssymb}
\usepackage{graphicx}
\newcommand{\Floor}[1]{\lfloor {#1} \rfloor}
\theoremstyle{definition}
\newtheorem{thm}{Theorem}[section]
\newtheorem{Lem}{Lemma}[section]
\newtheorem{Def}{Definition}[section]
\newcommand{\SL}{\text{SL}_2(\mathbb{Z})}
\newcommand{\HH}{\mathbb{H}}
\newcommand{\bb}[1]{\mathbb{{#1}}}
\title{Duke's Theorem and Continued Fractions}
\author{John Mangual}
\begin{document}
\maketitle
\abstract{  For uniformly chosen random $\alpha \in [0,1]$, it is known the probability the $n^{\rm th}$ digit of the continued-fraction expansion, $[\alpha]_n$  converges to the Gauss-Kuzmin distribution $\mathbb{P}([\alpha]_n = k) \approx \log_2 (1 + 1/ k(k+2))$ as $n \to \infty$.  In this paper, we show the continued fraction digits of $\sqrt{d}$, which are eventually periodic, also converge to the Gauss-Kuzmin distribution as $d \to \infty$ with bounded class number, $h(d)$.  The proof uses properties of the geodesic flow in the unit tangent bundle of the modular surface, $T^1(\text{SL}_2 \mathbb{Z}\backslash \mathbb{H})$.}
\section{Continued Fractions...}
\begin{paragraph}{} For any $\alpha \in [0,1]$ we can define the continued fraction expansion in $\mathbb{Z}$ by repeating a two step algorithm.  First $a_0 = \alpha$ and $b_0 = \Floor{\alpha}$.  Now we simply repeat:
\begin{equation}
a_{k+1} = \{1/a_k\} \hspace{2cm} b_{k+1} = \Floor{a_k} \end{equation}
The end result is that $\alpha$ can be encoded as a sequence of integers: $[b_0, b_1, b_2, \dots ]$. \end{paragraph}
\begin{paragraph}{} If $\alpha$ is rational then we get a finite continued fraction.  What if $\alpha$ is the square root of a irrational number?  Then we get an eventually repeating sequences of numbers $b_k$.  For example, the sequence for $\sqrt{7}$ is $[2,1,1,1,4,1,1,1,4, \dots]$ where the $[1,1,1,4]$ motif repeats forever.  How can we get a purely periodic sequence? A theorem by Galois says:
\end{paragraph}
\begin{thm} [5] A quadratic number $\alpha$ has a purely periodic continued fraction expansion if and only if $\alpha > 1$ and $-1 < \alpha' < 0$ where $\alpha, \alpha'$ have the same quadratic equation.  
\end{thm}
\begin{paragraph}{}Now let's ask about statistics of these continued fractions.  How often does the number $5$ appear in a generic continued fraction? This answer for a random $\alpha \in [0,1]$ chosen uniformly was found by Kuzmin in 1928.  He showed:  \end{paragraph}
\begin{thm}[5] There exist positive constants $A,B$ such that 
$$ \left| A_n(k) - \log_2\left( 1 + \frac{1}{k(k+2)}\right) \right| \leq \frac{A}{k(k+1)}e^{-B\sqrt{n-1}}$$ 
Where $A_n(k) = |\{ x \in [0,1]: b_n(x) = k \}|$.  \end{thm}
\begin{paragraph}{} In this paper we look at how the statitics of the continued fraction digits of $\sqrt{d}$ for $d > 0$ behave as $d \to \infty$.  In fact we have to be more specific and restrict ourselves to the case of bounded class number, so $h(d)$ is less than some constant.  Also, since our sequence $b_k(\sqrt{d})$ is deterministic, we need to define the statistics we'll be looking at:
$$ c(\alpha, k) = \lim_{T \to \infty} \frac{ \# \{0\leq i < T: b_i(\alpha) = k  \} }{T}$$
We claim that these statistics approach the limit above, i.e.
\end{paragraph}
\begin{thm} As $d \to \infty$ with $h(d)$ bounded:
$$\lim_{d \to \infty} c(\sqrt{d},k) \to \log_2 \left( 1 + \frac{1}{k(k+2)} \right) $$
\end{thm}
\begin{paragraph}{}
To prove this we're going, as $d \to \infty$ and $h(d)=1$, the orbits of $\sqrt{d}$ under the map $T: x \mapsto \{ 1/x\}$ approach the Gauss-Kuzmin on $[0,1]$.  We can rephrase Theorem 1.4 in this new language: \end{paragraph}
\begin{thm} Let $x_0 = \{\sqrt{d}\}$, $h(d) = 1$, $T:x \mapsto \{1/x\}$ be the Gauss map and $f:[0,1] \to \mathbb{R}$ be continuous:
$$ \lim_{d \to \infty}  \lim_{N\to \infty} \frac{1}{N}\sum_{k = 0}^{N-1}f(T^k(x_0))  = \int_0^1 \frac{f(x)}{\ln 2}\cdot \frac{dx}{1+x}$$
where $h(d) = 1$ as $d$ goes to infinity.\end{thm}
\begin{paragraph}{}  To prove this we need to change settings and examine geodesics in the upper half plane.\end{paragraph}

\section{... and the Geodesic Flow}
\begin{paragraph}{}  Let's switch contexts to $\HH = \{x + iy: y > 0\}$ as a differentiable manifold with 
the Poincar\'{e} metric: 
$$ ds^2 = \frac{dx^2 + dy^2}{y^2}$$
Much of this exposition follows [2], Chapters 3 and 13.\end{paragraph}
\begin{paragraph}{}
The geodesics in this metric are (Euclidean) semi-circles with diameters along the real line.  Thus for any unit tangent vector  $(z, e^{i\theta}) \in T^1(\HH)$ there is a unique oriented geodesic which goes through $z$ and whose tangent at $z$ points in the direction $e^{i\theta}$.  \end{paragraph}

\begin{paragraph}{} The group $\textrm{SL}(2,\mathbb{R})$ acts on $\HH$ by fractional linear transformations:
$$g = \left( \begin{array}{cc} a & b \\ c & d \end{array}\right): z \mapsto \frac{az+ b}{cz+d}$$
Then  points in the complex plane are identified as points in the projective complex line $\mathbb{P}^1(\mathbb{C}$:
$$\left(\begin{array}{c} z \\ 1 \end{array} \right) = \left( \begin{array}{c} \frac{az + b}{cz+ d} \\ 1 \end{array}\right)$$
This is simply the usual matrix action of $\mathrm{PSL}(2, \mathbb{R})$.
\end{paragraph}
\begin{paragraph}{} We can also consider the quotient group of $\mathbb{H}$ under the action of $\SL$.
The quotient under this group action $\SL \backslash \HH$ has is the fundamental domain represented by the intersection of four sets $\{\text{Im}(z) > 0\}$, $\{|z| < 1\}$, $\{|z + 1| < 1\}$ and $\{|z - 1| < 1\}$. 
\end{paragraph}
\begin{paragraph}{} The tangent space of $\mathbb{H}$ is simply $\HH \times \mathbb{C}$.  Any element $g \in \mathrm{SL}(2, \mathbb{R})$ can act on the tangent bundle by:
$$Dg(z,v)  = (g(z), g'(z) v)= \left( \frac{az+ b}{cz+d}, \frac{v}{(cz+d)^2}\right)$$ It turns out this action is simply transitive and therefore
$$T^1(\mathbb{H}) \simeq \mathrm{PSL}(2, \mathbb{R})$$
where $T^1(\mathbb{H})$ is the unit tangent bundle.
\end{paragraph}
\begin{paragraph}{}Now any two points, $z_1, z_2 \in \mathbb{H}$ determine a unique geodesic - the unique Euclidean semi-circle passing through $z_1$ and $z_2$.  We can consider a map, $\mathcal{G}_t$, which flows a tangent vector along its geodesic exactly $t$ units of arc length.  This is the {\it geodesic flow} from $z_1$ to $z_2$.  You need to specify both a starting point an a direction, an element of $S^1$, to get a unique geodesic.  Equivalently we can describe the geodesic as right-multiplication by the elements:

$$a_t = \left( \begin{array}{cc} e^{-t/2} & 0 \\ 0 & e^{t/2}  \end{array} \right)$$

This will be the basis for defining the continued fraction map in terms of the geodesic flow.  Furthermore we can define the geodesic flow restricted to $T^1(\SL)$.
\end{paragraph}

\begin{paragraph}{} In light of Theorem 1.1, we should consider geodesics whose endpoints $\alpha, \alpha' \in \bb{Q}[\sqrt{d}]$ for some $d > 0$ and such that $\alpha > 1$ and $-1 < \alpha' < 0$.  These curves necessarily 
cut the set $\{ yi : 0 < y < 1\}$ transversely.  In fact, we can identify these geodesics either by their endpoints or by the tangent vector at which they cut $[0,i]$.  
\end{paragraph}

\begin{paragraph}{}The elements of $B = T^1(\{ yi : 0 < y < 1\})$ parameterize these geodesics by the angle at which they cut the line segment $[0,i]$. If we only consider the first case, the set is called $B^+$ and in the second case it is called $B^-$.   
\begin{eqnarray*}
B^+ &=& \left\{ (yi, e^{i\theta}): 0 < y < 1, - \frac{\pi}{2} < \theta < 0 \right\} \\
B^- &=& \left\{ (yi, e^{i\theta}): 0 < y < 1, \pi < \theta < \frac{3\pi}{2} \right\}
\end{eqnarray*}
These correspond to purely periodic continued fractions and therefore to closed geodesics in the Riemann surface $\SL \backslash \bb{H}$.
\end{paragraph}
\begin{Def} Consider a geodesic, $\gamma$ in $\mathbb{H}$ with endpoints $\alpha, \alpha'$ for which $\alpha < -1$ and $1 > \alpha' > 0$. Define {\it natural coordinates} by 
$(y,z)$ where $y = \alpha$ and $z = \frac{1}{\alpha + \alpha'}$.
\end{Def}
\begin{Def} Define $T: B \to B$ in terms of the geodesic flow $\mathcal{G}_t$ by:
$$ T[(z, e^{i\theta})] = \mathcal{G}_{t_0} (z, e^{i\theta}) \text{\quad{}where\quad} t_0 = \inf\{ t > 0: \mathcal{G}_t(z, e^{i\theta}) \in \SL(B)\}$$ 
This $t_0$ may not be finite but whenever it is finite this map is well-defined.  This is known as the {\it return time} map for the cross section $B$.
\end{Def}
\begin{Lem}Let $x = (\mathrm{i}b, e^{i\theta} ) \in B_+$ have natural coordinates $(y,z)$. 
Then $T(x) \in \SL(B)$ if $T$ is defined on $x$.  Moreover $T(x) \in \SL(x')$, where $x'$ has natural coordinates
$$ \overline{T}(y,z) = \left( \left\{ \frac{1}{y}\right\}, y(1 - yz)\right)$$ A similar property holds for $B_-$. \end{Lem}
\begin{paragraph}{} The first step in the proof of Theorem 1.3 is to consider closed geodesics in $\SL \backslash \HH$ of discrimant $d$.  They correspond to purely periodic continued fractions and to elements of $B$.  Their natural coordinates lie in the ring $\bb{Q}[\sqrt{d}]$.  Thus, we can prove these geodesics become equidistributed in $B$ as $d \to \infty$ and $h(d) \leq M$, then the first natural coordinate follows the Gauss-Kuzmin distribution in $[0,1]$.  In other words, the orbit under $T: x \mapsto \{1/x\}$ of any element of $\bb{Q}[\sqrt{d}]$ becomes Gauss-Kuzmin
in $[0,1]$, asymptotically
\end{paragraph}

\section{Duke's Theorem}
\begin{paragraph}{} We to define some special sets of geodesics:
\end{paragraph}
\begin{Def} Let $d < 0$ and let $(x_d, y_d)$ be the fundamental solution to $x^2 - dy^2 = 4$.  Define $\Gamma_d$ as the set of geodesics in 
$\mathrm{PSL}_2(\mathbb{Z}) \backslash \mathbb{H}$ of length $d$ induced by quadratic forms $q(x,y) = ax^2 + bxy + cy^2$ with $b^2 - 4ac = d$. 
\end{Def}
\begin{paragraph}{}  We can consider closed geodesics on $\SL \backslash \HH$ 
\end{paragraph}
\begin{thm}[1] Suppose $d$ is a fundamental discriminant.  Then for some $\delta > 0$ depending only on $\Omega$ 
\begin{equation} \frac{ \sum_{C \in \Lambda_d} |C \cap \Omega|}{\sum_{C \in \Lambda_d} |C|} = \mu(\Omega) + \mathcal{O}(d^{-\delta})\end{equation} as $d \to \infty$ whee $|C|$ is the non-Euclidean length of $C$ and the $\mathcal{O}$ constant depends only on $\delta$ and $\Omega$. \end{thm}
\begin{paragraph}{} The idea that geodesic orbits in homogeneous spaces become equidistributed can be extended to the tangent bundle.  In this case, we consider the geodesic flow in $T^1\SL$.  This is useful because the set we wish to consider $B, B^+$ and $B^-$, who live in the unit tangent bundle and not the underlying space.  Fortunately for us, Duke's theorem extends to this case well.
\end{paragraph}
\begin{Def} Define $\Lambda_d = \{ \gamma_{[q]}: [q] \in \mathrm{PSL}_2(\mathbb{Z})\backslash Q_d(\mathbb{Z})\}$,  the geodesics associated
with the set of binary quadratic forms modulo $\mathrm{PSL}_2$ equivalence.  If $q(x,y) = ax^2 + bxy + cy^2$,  the geodesic $\gamma_{[q]}$ has endpoints defined by $q(x,1) = 0$.
Then project this geodesic onto the modular surface $\mathrm{PSL}_2(\mathbb{Z})\backslash\HH$.
\end{Def}
\begin{thm} [4] As $d \to \infty$, $d \equiv 0, 1 \mod 4$, $d$ not a perfect square, the set $\Gamma_d$ becomes equidistributed on the unit tangent bundle, $T^1\SL$,  with respect to the volume measure $d\mu_L = \frac{3}{\pi}\frac{dx dy}{y^2}\frac{d\theta}{2\pi}$. 
$$ \frac{ \sum_{C \in \Gamma_d} |C \cap \Omega|}{\sum_{C \in \Gamma_d} |C|} = \mu_L(\Omega) + \mathcal{O}(d^{-\delta})$$\end{thm}

\begin{paragraph}{}Now we are ready to prove our equidistribution result.  Theorem 1.4 can be rephrased in dynamical into ergodic theory language.   Proving it for $f(x) = \chi_I(x)$ for some interval $I \in [0,1]$, we can prove it for any continuous function $f(x)$,

\end{paragraph}
\begin{thm} Let $T:  [0,1]\to [0,1]$ be the map defined in Theorem 2.1.  Then the orbit $\{T^n(x): n \in \mathbb{N}\}$ is distributed as the Gauss-Kuzmin distribution as $d$ goes to infinity.  Specifically, for any interval $I \subseteq [0,1]$:
$$\lim_{d \to \infty} \frac{ \#\{0 \leq n < l(d):  T^n(B) \in I \}}{l(d)} = \frac{1}{\ln 2} \int_I\frac{ dx }{1 + x   }$$
Where $l(d)$ is the period of the continued fraction corresponding to $\sqrt{d}$.
\end{thm}
\begin{proof} Let $I \subseteq [0,1]$ be an interval and $(y,z)$ be the natural coordinates in $B$.
$$B_{I, \epsilon} = \{ (y,z) \in B: y \in I\} \times [-\epsilon/2, \epsilon/2]$$
This set can be embedded in $T^1(\SL \backslash \HH)$ as:
$$ \{ \mathcal{G}_t( B) : t \in [-\epsilon/2, \epsilon/2]\}$$
The Haar measure in the natural coodinates is Lebesgue in all three variables.

By Duke's theorem, there is $\delta > 0$ such that 
\begin{equation} \frac{\sum_{C \in \Gamma_D} | C \cap B_{I, \epsilon} | } {\sum_{C \in \Gamma_D} |C|}
= \mu(B_{I, \epsilon}) + \mathcal{O}(d^{-\delta}) 
 \end{equation}
Note that this equation is true even if $I = [0,1]$.  It therefore follows that:
\begin{equation} \frac{\sum_{C \in \Gamma_D} | C \cap B_{I, \epsilon}|}  {\sum_{C \in \Gamma_D} | C \cap B_{[0,1], \epsilon}|}= \frac{\mu(B_{I, \epsilon})}{ \mu(B_{[0,1], \epsilon} )} + \mathcal{O}(d^{-\delta}) 
 \end{equation}
 Because of how we defined $B_{I, \epsilon}$, the total length is simply the number of times the geodesics $\Gamma_d$
 cut $I$ times $\epsilon$:
 \begin{eqnarray*} \sum_{C \in \Gamma_D} | C \cap B_{I, \epsilon}| &=& \epsilon \cdot \sum_{C \in \Gamma_D} \# \{ C \cap I \} \\
&=& \epsilon \cdot   \# \{ 0 \leq n < l(d) : T^n(x_0) \in I \} \end{eqnarray*}
Here $l(d)$ denotes the {\it period} of the continued fraction with respect to $l(d)$. In the case $I = [0,1] $ 
the last line is just $\epsilon \cdot l(d)$ so the left hand side of (4) is really just counting measure:
$$ \frac{\# \{ 0 \leq n < l(d) : T^n(x_0) \in I \}   } { l(d)  }$$
For the right hand side of (4) let's find the measure of $B_{I, \epsilon}$:
$$\mu(B_{I, \epsilon}) = \epsilon \int_I \int_0^{1 + y}  \frac{dy dz}{\ln 2} =  \frac{\epsilon}{\ln 2}\int_I \frac{dy}{1 + y}$$
Therefore equation (4) should read:
$$ \lim_{d \to \infty} \frac{\# \{ 0 \leq n < l(d) : T^n(x_0) \in I \}   } { l(d)  }
= \frac{1}{\ln 2}\int_I \frac{dy}{1 + y}$$
\end{proof}

\begin{paragraph}{} That concludes the proof that these continued fractions of quadratic irrationalities
follow the Gauss-Kuzmin distribution as the discriminant tends to infinity.  Next we show this results describes truly generic behavior.
\end{paragraph}
\section{Bounded Class Number}
\begin{paragraph}{}We can relax the condition $h(d) = 1$ in Theorem 1.4: \end{paragraph}
 \begin{thm}    Let $x_0 = \{\sqrt{d}\}$, $h(d) = 1$, $T:x \mapsto \{1/x\}$ be the Gauss map and $f:[0,1] \to \mathbb{R}$ be continuous:
$$ \lim_{d \to \infty} \lim_{n\to \infty} \sum_{k = 0}^{n-1}f(T^k(x_0))  = \int_0^1 \frac{f(x)}{\ln 2}\cdot \frac{dx}{1+x}$$
where there exists $M > 0$ such that
$h(d) < M$. \end{thm}
\begin{proof} We simply examine Theorem 3.3 more closely.  In this case, $\Gamma_d$ has more than one element, i.e. there are several geodesics with the same discrimant $d$. Duke's theorem says the sum
$$ \frac{ \sum_{C \in \Lambda_d} |C \cap \Omega|}{\sum_{C \in \Lambda_d} |C|} = 
\frac{ \sum_{C \in \Lambda_d} |C \cap \Omega|}{|\Lambda_d| \cdot |L_d|}
= \mu(\Omega) + \mathcal{O}(d^{-\delta})$$
is the Lebesgue measure.  Here $L_d$ is the length of a geodesic of discriminant $d$.  Following the principle outlined in [6] (Section 1.3.5 (1)), since the Lebesgue measure is an extreme point in the convex space of geodesic-flow invariant measures on $T^1(\SL)$, since the Gauss map $T$ is ergodic and since their sum is a Lebesgue measure, 
each term in the sum must also approach Lebesgue measure.  
\end{proof}
\begin{paragraph}{}  Then continued fractions in $\bb{Q}[\sqrt{d}]$ tend towards the Gauss Kuzmin distribution even in the case bounded class number.  This behavior is truly generic.
\end{paragraph}

\end{document}